\documentclass[oneside,12pt]{article}
\usepackage{graphicx}
\begin{document}
\title{Heegaard Splittings of $\partial$-reducible 3-manifolds}

\author{Jiming Ma$^{*}$ and Ruifeng Qiu\footnote{Both authors are supported by a grant (No.10171038) of NSFC}}
\maketitle \vskip 3mm

{\bf Abstract.} In this paper, we shall prove that any Heegaard
splitting of a $\partial$-reducible 3-manifold $M$, say $M=W\cup
V$, can be obtained by doing connected sums, boundary connected
sums and self-boundary connected sums from Heegaard splittings of
$n$ manifolds $M_{1},\ldots, M_{n}$ where $M_{i}$ is either a
solid torus or a $\partial$-irreducible manifold. Furthermore,
$W\cup V$ is stabilized if and only if one of the factors is
stabilized.\vskip 3mm

{\bf Keywords.} Connected sum, boundary connected sum,
self-boundary connected sum.\vskip 3mm

AMS Subject Classification: 57N10, 57M50. \vskip 3mm

{\bf \S 1 Introduction}\vskip 3mm

Let $M$ be a compact 3-manifold with boundary such that each
component of $\partial M$ is not a 2-sphere. If there is a
2-sphere in $M$ which does not bound any 3-ball, then we say $M$
is reducible; otherwise, $M$ is irreducible. If there is an
essential disk $D$ in $M$, then we say $M$ is
$\partial$-reducible; otherwise, $M$ is $\partial$-irreducible.

Let $M$ be a compact 3-manifold. If there is a closed surface $S$
which separates $M$ into two compression bodies $W$ and $V$ with
$\partial_{+} W=\partial_{+} V=S$, then we say $M$ has a Heegaard
splittings, denoted by $M=W\cup_{S}V$ or $M=W\cup V$. In this
case, $S$ is called a Heegaard surface of $M$. A Heegaard
splitting $M=W\cup_{S} V$ is said to be reducible if there are two
essential disks $D_{1}\subset W$ and $D_{2}\subset V$ such that
$\partial D_{1}=\partial D_{2}$. A Heegaard splitting $M=W\cup_{S}
V$ is said to be $\partial$-reducible if there is an essential
disk $D$ in $M$ such that $D$ intersects $S$ in only one essential
simple closed curve. A Heegaard splitting $M=W\cup_{S} V$ is said
to be weakly reducible if there are two essential disks
$D_{1}\subset W$ and $D_{2}\subset V$ such that $\partial
D_{1}\cap \partial D_{2}=\emptyset$. A Heegaard splitting
$M=W\cup_{S} V$ is said to be stabilized if there are two properly
embedded disks $D_{1}\subset W$ and $D_{2}\subset V$ such that
$D_{1}$ intersects $D_{2}$ in only one point. It is easy to see
that if $M=W\cup_{S} V$ is stabilized and $g(S)\geq 2$ then it is
reducible.

Now there are some results on reducibilities of Heegaard
splittings. For example, Haken proved that any Heegaard splitting
of a reducible 3-manifold is reducible; Casson and Gordon gave a
disk version of Haken's lemma, that say, any Heegaard splitting of
a $\partial$-reducible 3-manifold is $\partial$-reducible, they
also show that if $M$ has a weakly reducible Heegaard splitting
$W\cup V$ then either $W\cup V$ is reducible or $M$ contains an
essential closed surface of genus at least one; Ruifeng Qiu
recently proved Gordon's conjecture on stabilizations of reducible
Heegaard splittings, that say, the connected sum of two Heegaard
splittings is stabilized if and only if one of the two factors is
stabilized.

In this paper, we shall consider  Heegaard splittings of
$\partial$-reducible manifolds. The main result is the following:

{\bf Theorem 1.} \  Any Heegaard splitting of a
$\partial$-reducible manifold $M$, say $M=W\cup V$, can be
obtained by doing connected sums, boundary connected sums and
self-boundary connected sums from Heegaard splittings of $n$
manifolds $M_{1},\ldots, M_{n}$, where $M_{i}$ is either a solid
torus or an irreducible, $\partial$-irreducible manifold.
Furthermore, $W\cup V$ is stabilized if and only if one of the
factors is stabilized.

{\bf Remark.} \ If $M=W\cup V$ can be obtained by doing connected
sums, boundary connected sums and self-boundary connected sums
from Heegaard splittings of $l$ manifolds $X_{1},\ldots, X_{l}$
where $X_{j}$ is either a solid torus or an irreducible,
$\partial$-irreducible manifold, then $n=l$ and
$\bigl\{M_{1},\ldots,
M_{n}\bigr\}=\bigl\{X_{1},\ldots,X_{n}\bigr\}$. We omit the
proof.\vskip 3mm

A Heegaard splitting of a handlebody $H$, say  $W\cup V$, is said
to be trivial if $W$ is homeomorphic to $\partial H\times I$ and
$V$ is homeomorphic to $H$.

As an application, we shall give a new proof to
Scharlemann-Thomston's result:

{\bf Corollary 2([ST]).} \ Any unstabilized Heegaard splitting of
a handlebody is trivial.\vskip 3mm

{\bf \S 2 Prelimary}\vskip 3mm

{\bf  Connected sums of Heegaard splittings}

Now let $M=W\cup V$ be a reducible Heegaard splitting. Then there
is a 2-sphere $P$ such that $B_{W}=P\cap W$ is an essential disk
in $W$ and $B_{V}=P\cap V$ is an essential disk in $V$.  Suppose
that $P$ separates $M$ into $M_{1}^{*}$ and $M_{2}^{*}$. Then $B_{
W}$ separates $W$ into $ W_{1}^{*}$ and $W_{2}^{*}$, $B_{V}$
separates $V$ into $V_{1}$ and $V_{2}$. We may assume that
$W_{1}^{*}, V_{1}\subset$$ M_{1}^{*}$ and $W_{2}^{*}, V_{2}\subset
M_{2}^{*}$. Let $M_{1}=M_{1}^{*}\cup_{P} H^{3}_{1}$ and
$M_{2}=M_{2}^{*}\cup_{P} H^{3}_{2}$ where $H^{3}_{1}$ and
$H^{3}_{2}$ are two 3-balls. Then $M$ is called the connected sum
of $M_{1}$ and $M_{2}$, denoted by $M=M_{1}\sharp M_{2}$. Let
$W_{1}=W_{1}^{*}$$\cup H^{3}_{1}$, and $W_{2}=W_{2}^{*}$$\cup
H^{3}_{2}$. Then $W_{1}$ and $W_{2}$ are two compression bodies
such that $\partial_{+} V_{1}=\partial_{+} W_{1}$ and
$\partial_{+} V_{2}=\partial_{+} W_{2}$. Hence $M_{1}=W_{1}\cup
V_{1}$ is a Heegaard splitting of $M_{1}$ and $M_{2}=W_{2}\cup
V_{2}$ is a Heegaard splitting of $M_{2}$. In this case, $W\cup V$
is called the connected sum of $ W_{1}\cup V_{1}$ and $W_{2}\cup
V_{2}$.\vskip 3mm

{\bf Boundary connected sums of Heegaard splittings}

Let $M$ be a compact orientable $\partial$-reducible 3-manifold,
and $D$ be an essential disk in $M$. Suppose that $D$ is
separating in $M$. Then  $\partial D$ is also separating in
$\partial M$. Now each component of $M-D\times (0,1)$ is a
3-manifold with boundary, denoted by  $M_i$. Without loss of
generality, we may assume that $D\times\bigl\{0\bigr\}\subset
\partial M_{1}$ and $D\times\bigl\{1\bigr\}\subset \partial M_{2}$.
In this case, we say $M$ is the boundary connected sum of $M_{1}$
and $M_{2}$, denoted by $M= M_1 \cup_{D}M_2$.

Suppose that $M_i=W^i \cup V^i$ is a Heegaard of $M_i$ such that
$D\times \bigl\{0\bigr\}\subset
\partial_{-}V^1$ and $D\times
\bigl\{1\bigr\}\subset
\partial_{-}V^2$. Then there are unknotted, properly embedded arcs
$\alpha_i$ in $V^i$ and  $\beta$ in $D\times [0,1]$  such that
$\partial_{1} \alpha_{1}\subset \partial_{+} V^{1}$,
$\partial_{2}\alpha_{1}=\partial_{1}\beta$,
$\partial_{2}\beta=\partial_{1}\alpha_{2}$, $\partial_{2}
\alpha_{2}\subset\partial_{+} V^{2}$. Then
$\gamma=\alpha_{1}\cup\beta\cup\alpha_{2}$ is a properly embedded
arc in $V^{1}\cup D\times [0,1]\cup V^{2}$. Now let $N(\gamma)$ be
a regular neighborhood $\gamma$ in $V^{1}\cup D\times [0,1]\cup
V^{2}$ such that $N(\partial_1 \gamma)\subset
\partial_{+} V^{1}$, and $N(\partial_2 \gamma)\subset \partial_{+}
V^{2}$. It is easy to see that $W^1 \cup N(\gamma)\cup W^2)$,
denoted by $W$, is a compression body in $M$, and the closure of
$V^1\cup V^2-N(\gamma)$, denoted by $V$, is also a compression
body
 in $M$. Hence $W\cup V$ is a  Heegaard splitting of $M$. We say  $W\cup V$ is the
 boundary connected sum of the two Heegaard splittings $W^{1}\cup V^{1}$ and $W^{2}\cup
 V^{2}$.\vskip 3mm

{\bf Self-boundary connected sums of Heegaard splittings}

Let $M$ be a $\partial$-reducible 3-manifold, and $D$ be an
essential disk in $M$. Suppose that $D$ is non-separating in $M$,
but $\partial D$ is separating in $\partial M$. Now
$M^{'}=M-D\times (0,1)$ is a connected manifold such that
$\partial M^{'}$ contains at least two components $F_{1}$ and
$F_{2}$. We may assume that $D\times \bigl\{0\bigr\}\subset F_{1}$
and $D\times \bigl\{1\bigr\}\subset F_{2}$. In this case, we say
that $M$ is a self-boundary connected sum of $M^{'}$.
\begin{center}
\includegraphics[totalheight=6.5cm]{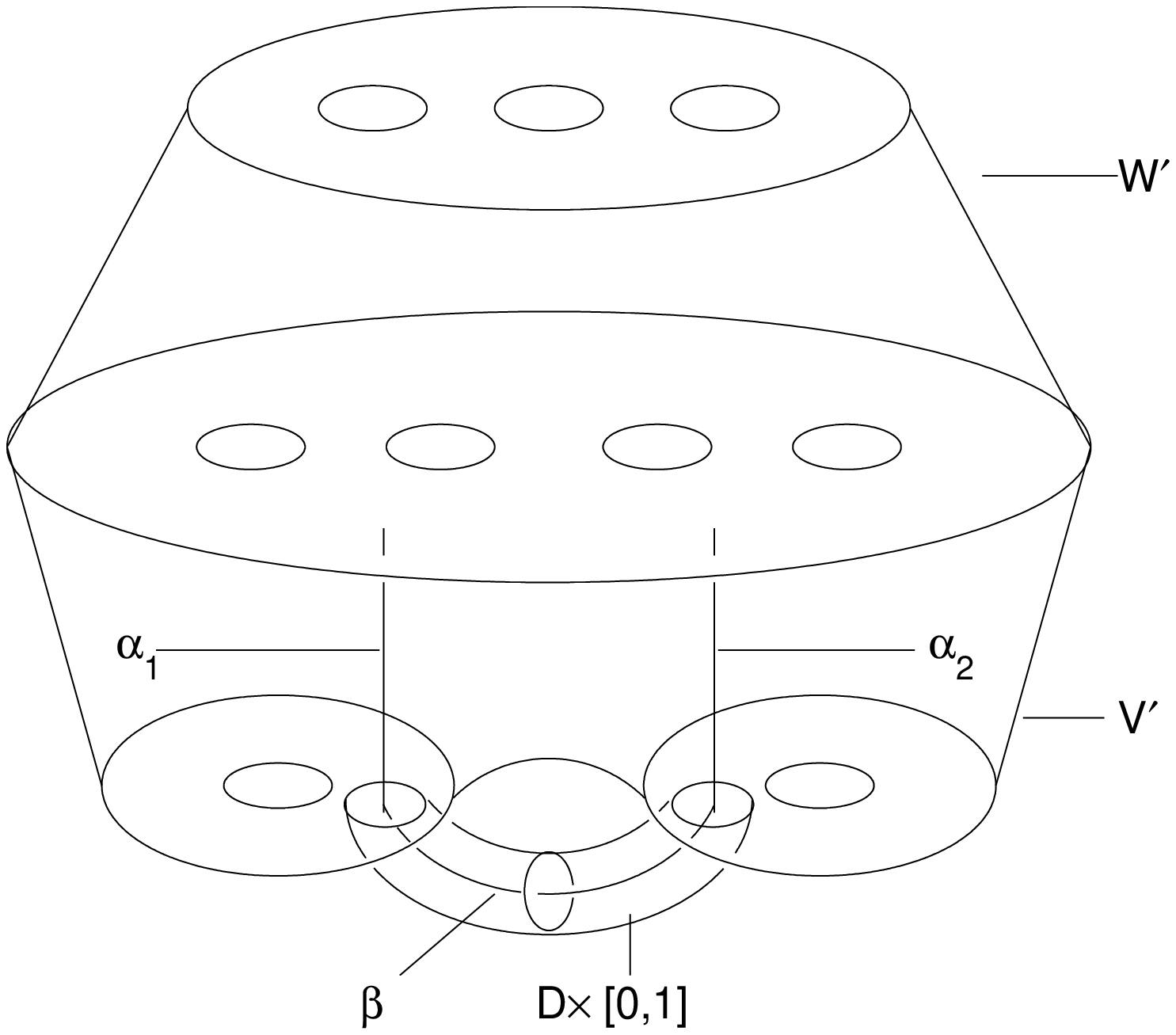}
\begin{center} Figure 1
\end{center}
\end{center}

\begin{center}
\includegraphics[totalheight=7cm]{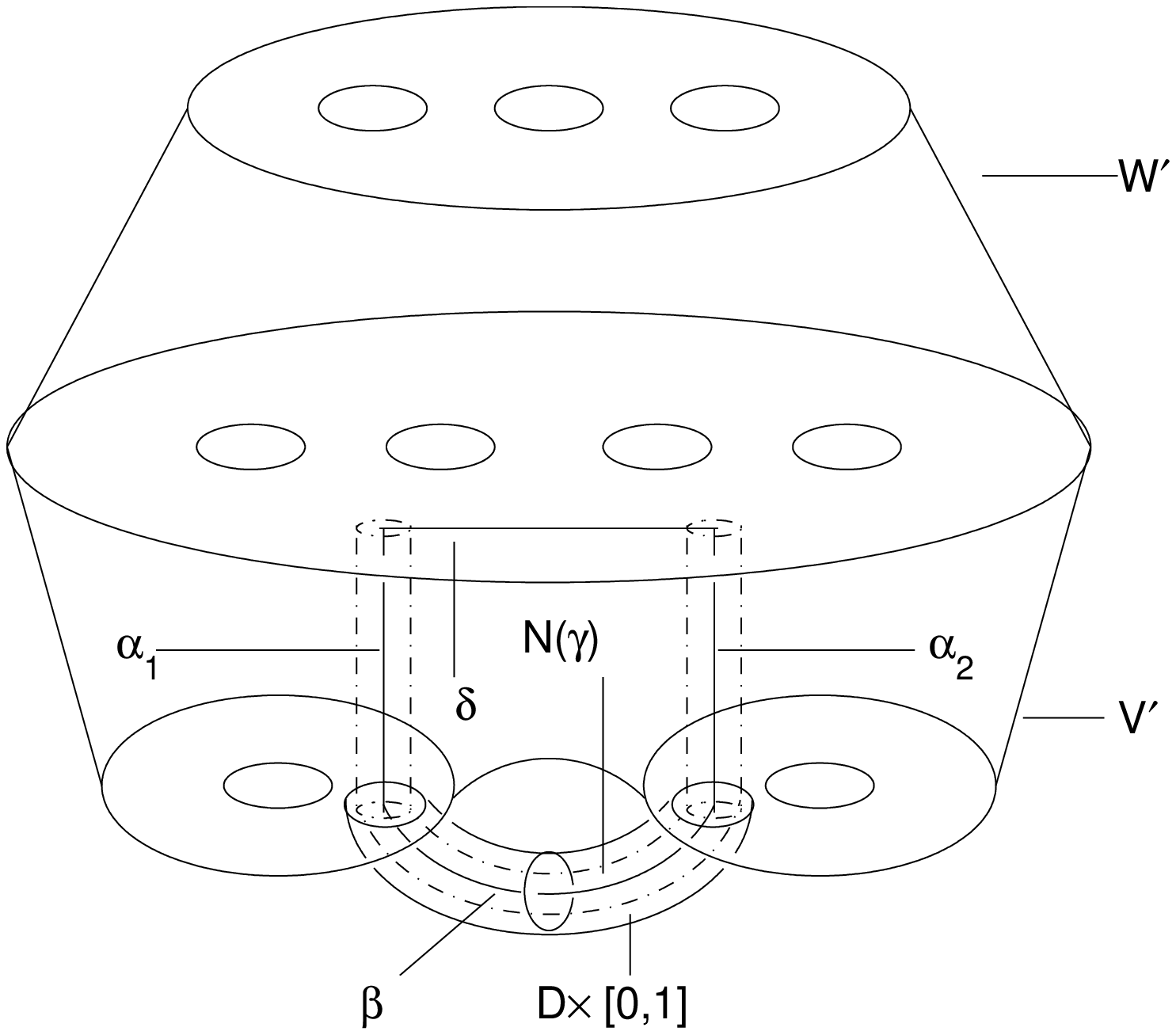}
\begin{center} Figure 2
\end{center}
\end{center}

let $M^{'}=W^{'}\cup V^{'}$ be a Heegaard splitting of $M^{'}$,
such that $F_1, F_2 \subset
\partial_{-} V^{'}$. Now suppose that $\alpha_{1}$, $\alpha_{2}$ are two unknotted
properly embedded arc in $V^{'}$ such that $\partial_{1}\alpha_{1}
$ and $\partial_{2}\alpha_{2}$ lie in $\partial_{+}V^{'}$, and
$\beta$ is a unknotted properly embedded arc in $D\times [0,1]$
such that $\partial_{2}\alpha_{1}=\partial_{1}\beta$ and
$\partial_{1}\alpha_{2}=\partial_{2}\beta$. Then
$\gamma=\alpha_{1}\cup\beta\cup\alpha_{2}$ is a properly embedded
arc in $V^{'}\cup D\times [0,1]$. Let $N(\gamma)$ be a regulalr
neighborhood of $\gamma$  in $V^{'}\cup D\times [0,1]$. It is easy
to see that $W=W^{'}\cup N(\gamma)$ is a compression body and the
closure of $V^{'}\cup D\times [0,1]-N(\gamma)$, denoted by $V$, is
also a compression body. Hence $M=W\cup V$ is a Heegaard splitting
of $M$. We say $W\cup V$ is a self-boundary connected sum of
$W^{'}\cup V^{'}$. See Figure 1 and Figure 2.

By definitions, if $W\cup_{S}V$ is the connected sum or the
boundary connected sum of $W^{1}\cup_{S_{1}}V^{1}$ and
$W^{2}\cup_{S_{2}}V^{2}$, then $g(S)=g(S_{1})+g(S_{2})$; if
$W\cup_{S} V$ is the self-boundary connected sum of
$W^{'}\cup_{S^{'}} V^{'}$, then $g(S)=g(S^{'})+1$.\vskip  3mm

{\bf \S 3 Proofs of Theorem 1}\vskip 3mm

{\bf Lemma 3.1 ([Q]).} \ The connected sum of two Heegaard
splittings is stabilized if and only if one of the two factors is
stabilized.

\begin{center}
\includegraphics[totalheight=7cm]{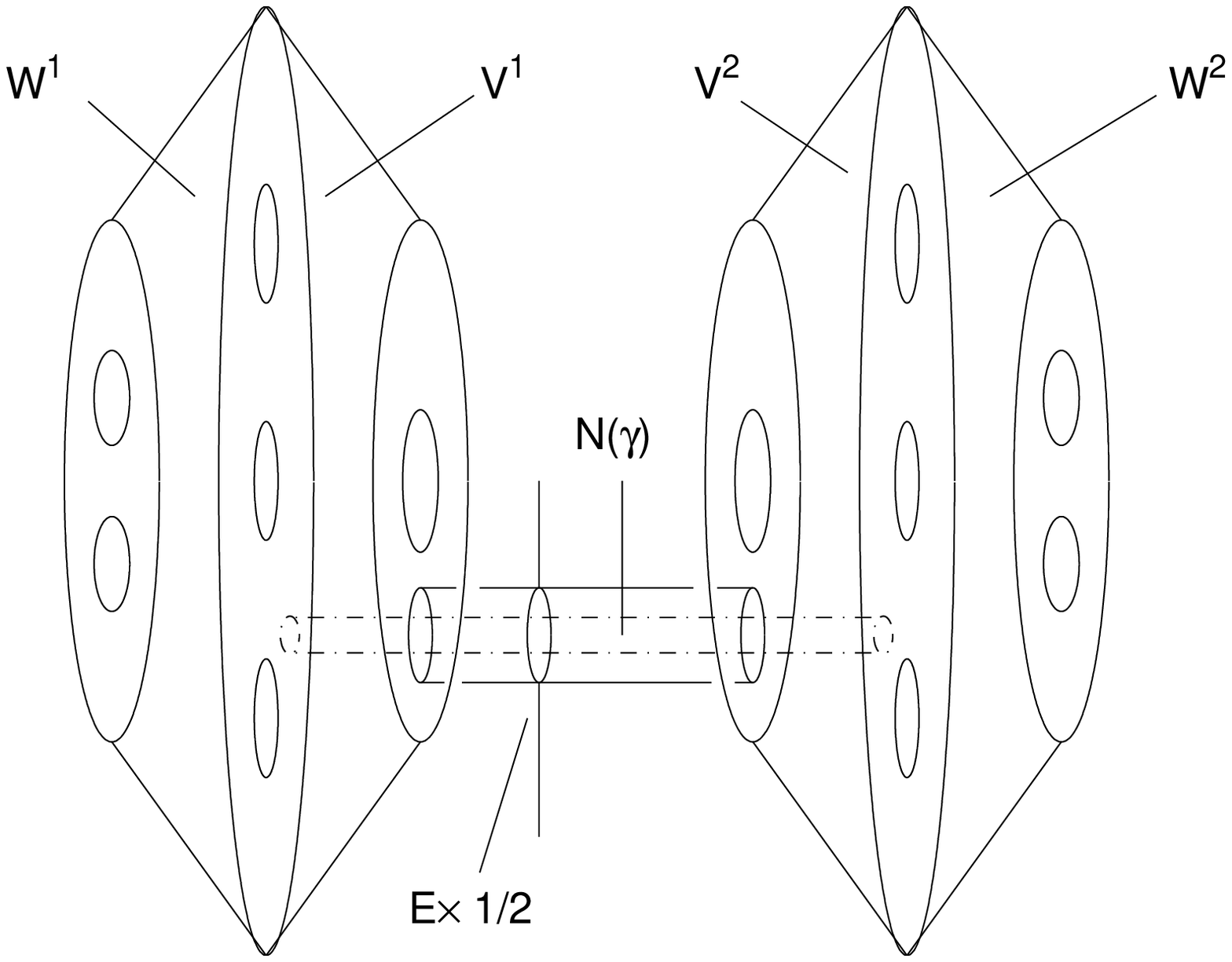}
\begin{center} Figure 3
\end{center}
\end{center}\vskip 3mm

{\bf Lemma 3.2.}  The boundary connected sum of two Heegaard
splittings is stabilized if and only if one of the two factors is
stabilized.

{\bf Proof.} Suppose that $M_1$ and $M_2$ are two 3-manifolds with
boundary, and $M=M_{1}\cup D\times [0,1]\cup M_{2}$ such that
$D\times\bigl\{0\bigr\}\subset
\partial M_{1}$ and $D\times\bigl\{1\bigr\}\subset \partial
M_{2}$. Suppose that $M_i=W^i \cup V^i$ is a  Heegaard splitting
of $M_i$ for $i=1,2$,  and $M=W\cup_{S} V$ is the boundary
connected sum of $W^{1}\cup V^{1}$ and $W^{2}\cup V^{2}$ as in
Figure 3.

Assume first that one of $W^{1}\cup V^{1}$ and $W^{2}\cup V^{2}$,
say $W^{1}\cup V^{1}$, is stabilized. Then there are two essential
disks $B_{W}\subset W^{1}$ and $B_{V}\subset V^{1}$  such that
$|B_{W}\cap B_{V}|=1$. It is easy to see that $B_{W}$ and $B_{V}$
can be chosen to be disjoint from $N(\alpha_{1})\subset
N(\gamma)$. Hence $W\cup_{S} V$ is also stabilized.

Assume now that $W^{1}\cup V^{1}$ and $W^{2}\cup V^{2}$ are
unstabilized. We attach a 2-handle $E\times I$ to $M$ such that
$\partial E\times I=\partial D\times I$, where $E$ is a disk. We
debote the resulting manifold  by $M^{*}$. Let $V^{*}=V\cup
E\times I$, then $V^{*}$ is a compression body. Hence $M^{*}=M\cup
E\times I$ has a Heegaard splitting $M^{*}=W\cup_S V^{*}$. Since
$D\times\bigl\{t\bigr\}$ intersects $S$ in an essential simple
closed curve which lies in $\partial N(\gamma)$,
$P=D\times\bigl\{1/2\bigr\}\cup E\times\bigl\{1/2\bigr\}$ is a
separating 2-sphere in $M$ which intersects  $S$ in an essential
simple closed curve.  By definition,  $M^*$ is just the connected
sum of $M_1$ and $M_2$ and $M^*=W\cup V^*$  is the connected sum
of $W^{1}\cup V^{1}$ and $W^{2}\cup V^{2}$. Since $W^{i}\cup
V^{i}$ is unstabilized, by Lemma 3.1, $W\cup V^{*}$ is
unstabilized. Hence $W\cup V$ is unstabilized. Q.E.D.\vskip 4mm

{\bf Lemma 3.3.} \ The self-boundary connected sum of
 a Heegaard splitting $M^{'}=W^{'}\cup V^{'}$ is stabilized if and only if
$M^{'}=W^{'}\cup V^{'}$ is stabilized.

{\bf Proof.} Suppose that $M^{'}=W^{'}\cup V^{'}$ is a  Heegaard
splitting of $M^{'}$ and $M=W\cup_{S} V$ is a self-boundary
connected sum of $M^{'}=W^{'}\cup V^{'}$ defined in Section 2.

Assume first that $W^{'}\cup V^{'}$  is stabilized. Then there are
two essential disks $B_{W}\subset W^{'}$ and $B_{V}\subset V^{'}$
such that $|B_{W}\cap B_{V}|=1$. It is easy to see that $B_{W}$
and $B_{V}$ can be chosen to be disjoint from $N(\alpha_{1}\cup
\alpha_{2})\subset N(\gamma)$. Hence $W\cup_{S} V$ is also
stabilized.

Assume now that $W^{'}\cup V^{'}$ is unstabilized. We attach a
2-handle $E\times I$ to $V$ such that $\partial E\times I=\partial
D\times I$, where $E$ is a disk. Let $M^{*}=M\cup E\times I$. Then
$P=D\times \bigl\{1/2\bigr\} \cup E\times \bigl\{1/2\bigr\}$ is a
non-separating 2-sphere in $M^{*}$. By definition, $P$ intersects
$S$ in an essential simple closed curve.

Let $\delta$ be an arc in $\partial_{+}W^{'}$, such that $\partial
\delta =\partial\gamma$ as in Figure 2, then $\gamma\cup \delta$
is a simple closed curve which intersects $P$ in one point. Let
$P\times I$ be a regular neighborhood of $P$ in $E\times I\cup
D\times I$. Let $a=\gamma\cup\delta-P\times (0,1)$ and $N(a)$ be a
regular neighborhood of $a$ in $M^{*}-P\times (0,1)$ such that

1) \ each  of $N(\delta)\cap W^{'}$ and $N(\delta)\cap V^{'}$ is a
half 3-ball,

2) \ $N(a)\cap N(\gamma)\subset N(a)$ where $N(\gamma)$ is defined
in Section 2.

Now let $V^{*}=V\cup E\times I$. Then $M^{*}=W\cup V^{*}$ is a
Heegaard splitting of $M^{*}$. Let $P^{*}=\partial N(a)\cup
P\times\bigl\{0,1\bigr\}-int(N(a)\cap(P\times\bigl\{0,1\bigr\}))$.
Then $P^{*}$ is a 2-sphere which intersects $S=\partial_{+}
W=\partial_{+} V^{*}$ in an essential simple closed curve. Now by
observations, $M^{*}=W\cup V^{*}$ is the connected sum of
$M^{'}=W^{'}\cup V^{'}$ and a genus one Heegaard splitting of
$S^{2}\times S^{1}$ along $P^{*}$.

Since $M^{'}=W^{'}\cup V^{'}$ is unstabilized, by Lemma 3.1,
$M^{*}=W\cup V^{*}$ is unstabilized. Note that $V^{*}=V\cup
E\times I$. Hence $M=W\cup V$ is also unstabilized. Q.E.D.\vskip
4mm

\begin{center}
\includegraphics[totalheight=7cm]{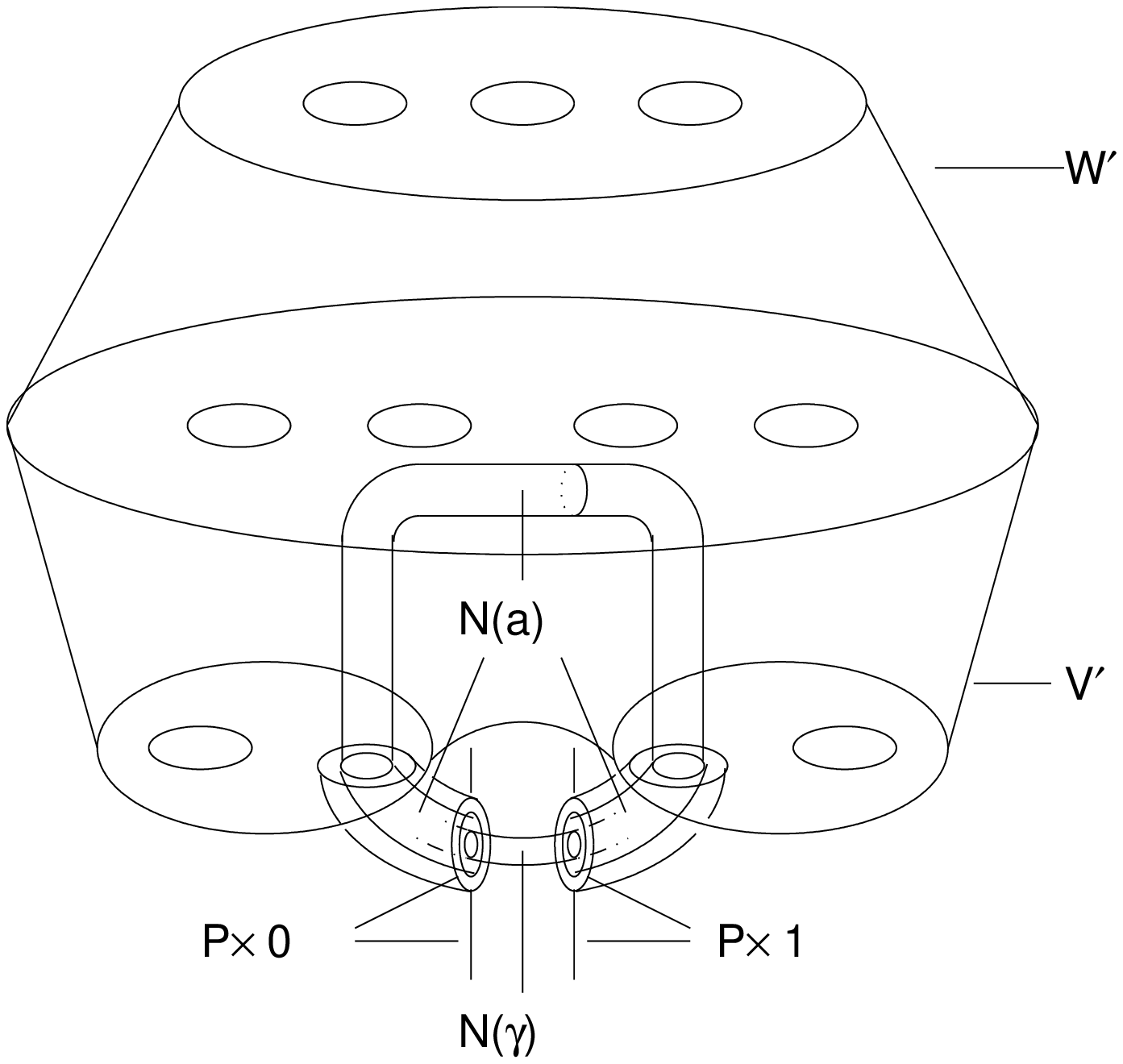}
\begin{center} Figure 4
\end{center}
\end{center}\vskip 3mm

Lemma 3.2 and Lemma 3.3 mean that  Heegaard genus is additive
under boundary connected sums and self-boundary connected sums
even if Heegaard genus is not minimal.\vskip 3mm

{\bf Lemma 3.4 ([CG]).} \ Any Heegaard splitting of a
$\partial$-reducible 3-manifold is $\partial$-reducible.\vskip 3mm

{\bf Theorem 1.}   Any Heegaard splitting of a
$\partial$-reducible manifold $M$, say $M=W\cup_{S} V$, can be
obtained by doing connected sums, boundary connected sums and
self-boundary connected sums from Heegaard splittings of $n$
manifolds $M_{1},\ldots, M_{n}$ where $M_{i}$ is either a solid
torus or a $\partial$-irreducible manifold. Furthermore,
$M=W\cup_{S} V$ is stabilized if and only if one of the factors is
stabilized.

{\bf Proof.} \ Suppose that $M=W\cup_{S} V$ is a Heegaard
splitting of a $\partial$-reducible 3-manifold. If the genus of
$M=W\cup_{S} V$ is one, then  $M$ is a solid torus and
$M=W\cup_{S} V$ is a trivial Heegaard splitting of $M$. So we may
assume that the genus of $M=W\cup_{S} V$ is at least two.

By Lemma 3.4, there is an essential disk $D$ such that $D$
intersects $S$ in an essential simple closed curve in $S$. We may
assume that $D\cap W$ is a disk and $D\cap V$ is an annulus. That
means that $\partial D\subset \partial_{-} V$. Now there are three
cases:

{\bf Case 1.} \ $D$ is separating in $M$.

Now $M-D\times (0,1)$ contains two components $M_{1}$ and $M_{2}$,
$D\cap W$ separates $W$ into two compression bodies $W_{1}$ and
$W_{2}$ and $D\cap V$ separates $V$ into two components $V_{1}$
and $V_{2}$. We assume that $W_{1},V_{1}\subset M_{1}$ and
$W_{2},V_{2}\subset M_{2}$. Let $N(D\cap W\times \bigl\{0\bigr\})$
be a regular neighborhood of $D\cap W\times \bigl\{0\bigr\}$ in
$W_{1}$ and $N(D\cap W\times \bigl\{1\bigr\})$ be a regular
neighborhood of $D\cap W\times \bigl\{1\bigr\}$ in $W_{2}$. Then
$V_{1}\cup N(D\cap W\times \bigl\{0\bigr\})$ and $V_{2}\cup
N(D\cap W\times \bigl\{1\bigr\})$ are two compression bodies.
Hence $(W_{1}-N(D\cap W\times\bigl\{0\bigr\}))\cup(V_{1}\cup
N(D\cap W\times \bigl\{0\bigr\})$ is a Heegaard splitting of
$M_{1}$ and $(W_{2}-N(D\cap W\times\bigl\{1\bigr\}))\cup(V_{2}\cup
N(D\cap W\times \bigl\{1\bigr\})$ is a Heegaard splitting of
$M_{2}$. By definition, $W\cup V$ is a boundary connected sum of
$(W_{1}-N(D\cap W\times\bigl\{0\bigr\}))\cup(V_{1}\cup N(D\cap
W\times \bigl\{0\bigr\})$ and $(W_{2}-N(D\cap
W\times\bigl\{1\bigr\}))\cup(V_{2}\cup N(D\cap W\times
\bigl\{1\bigr\})$.

By Lemma 2.2, $W\cup V$ is stabilized if and only if one of
$(W_{1}-N(D\cap W\times\bigl\{0\bigr\}))\cup(V_{1}\cup N(D\cap
W\times \bigl\{0\bigr\})$ and $(W_{2}-N(D\cap
W\times\bigl\{1\bigr\}))\cup(V_{2}\cup N(D\cap W\times
\bigl\{1\bigr\})$  is stabilized.\vskip 3mm

{\bf Case 2.} \ $D$ is non-separating in $M$, but $\partial D$ is
separating in $\partial M$.

{\bf Claim 1.} \ $D\cap W$ is non-separating in $W$.

{\bf Proof.} \ Suppose, otherwise, that $D\cap W$ is separating in
$W$. Then $\partial (D\cap W)$ is separating in $\partial_{+}
W=\partial_{+} V$. Since $V$ is a compression body, $D\cap V$ is
separating in $V$. Hence $D$ is separating in $M$, a
contradiction. Q.E.D. (Claim 1)

Now $M-D\times (0,1)$ is a manifold $M^{'}$.  Since $D\cap W$ is a
non-separating disk in $W$,  $W-(D\cap W) \times (0,1)$ is a
compression body, say $W^{*}$.  Let $N((D\cap W)\times
\bigl\{0\bigr\})$ be a regular neighborhood of $(D\cap W)\times
\bigl\{0\bigr\}$  and $N((D\cap W)\times \bigl\{1\bigr\})$ be a
regular neighborhood of $(D\cap W)\times \bigl\{1\bigr\}$ in
$W^{*}$. Then $(V-D\times (0,1))\cup N((D\cap W)\times
\bigl\{0\bigr\})\cup N((D\cap W)\times \bigl\{1\bigr\})$ is a
compression body, say $V^{'}$, in $M{'}$. Note that the closure of
$W^{*}-(N((D\cap W)\times \bigl\{0\bigr\})\cup N((D\cap W)\times
\bigl\{1\bigr\}))$, say $W^{'}$, is also a compression body. By
definition, $W\cup V$ is a self-boundary connected sum of
$W^{'}\cup V^{'}$.

By Lemma 3.3, $W\cup V$ is stabilized if and only if $W^{'}\cup
V^{'}$ is stabilized.

{\bf Case 3.} \ $D$ is non-separating in $M$, and $\partial D$ is
non-separating in $\partial M$.

{\bf Claim 2.} \ $\partial (D\cap W)$ is non-separating in
$S=\partial_{+} V=\partial_{+} W$.

{\bf Proof.} \ Suppose, otherwise, that $\partial (D\cap W)$ is
separating in $S$. Without loss of generality, we may assume that
$\partial_{-} V$ contains only one component. Let $V^{*}$ be the
manifold obtained by attaching a handlebody $H$ to $V$ along
$\partial_{-} V$ such that $\partial D$ bounds a disk $D^{*}$ in
$H$. Then $V^{*}$ is a handlebody and $(D\cap V)\cup D^{*}$ is a
disk in $V^{*}$. Since $\partial (D\cap W)$ is separating in $S$,
$(D\cap V)\cup D^{*}$ is separating in $V^{*}$, but $D^{*}$ is
non-separating in $H$, a contradiction. Q.E.D. (Claim 2)
\begin{center}
\includegraphics[totalheight=4.5cm]{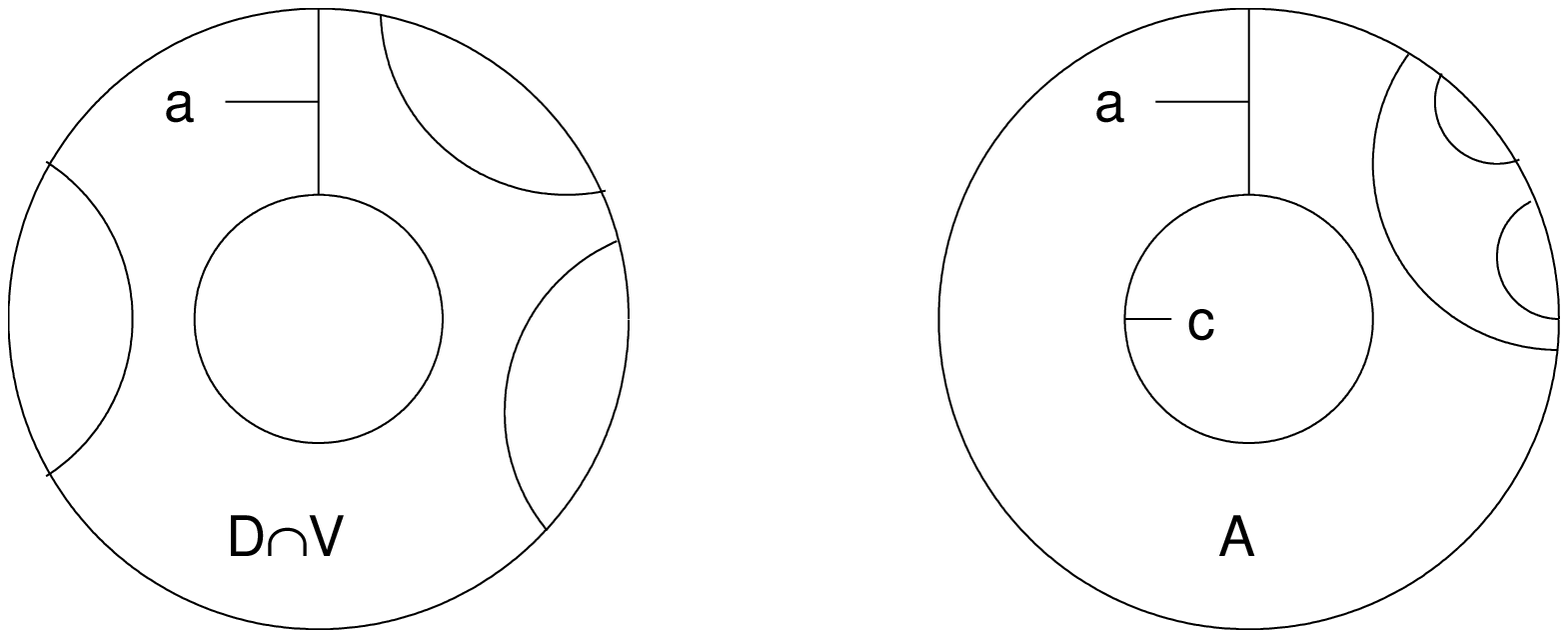}
\begin{center} Figure 5
\end{center}
\end{center}\vskip 3mm

{\bf Claim 3.} \ There is an annulus $A$ such that

1) \ one component of $A$ lies in $\partial_{+} V$ and the other
lies in $\partial_{-} V$, and

2) \ $A$ intersects the annulus $D\cap V$ in only one essential
arc.

{\bf Proof.} \ Now since $\partial_1 (D\cap V)$ in $\partial_- V$
is a non-separating curve, there is a curve in $\partial_- V$, say
$c$, such that $\mid\partial_1 (D\cap V)\cap c\mid=1 $. Then $c$
together with a simple closed curve in $\partial_+ V$ cobound an
annulus, say $A$. We may assume that $\mid A\cap (D\cap V)\mid$ is
minimal among all such annuli. Now we  prove $\mid A\cap (D\cap
V)\mid=1$.

Note that $A$ and $D\cap V$ are incompressible in $V$. Hence  $
A\cap (D\cap V)$ is a set of arcs. Since $c$ intersects
$\partial_{1}(D\cap V)$ in one point, there is only one arc, say
$a$, in $ A\cap(D\cap V)$ which is essential in both $A$ and
$D\cap V$.  See Figure 5.

Suppose that $\mid A\cap (D\cap V)\mid>1$. Let $b$ be an arc in
$A\cap(D\cap V)$ which is outermost in $D\cap V$, then it, with a
sub-arc of $\partial_2 (D\cap V)$, cobound a disk $E$ in $D\cap V$
such that $intE$ is disjoint from $A$. Now $b$, with a sub-arc of
$\partial_2 A$, cobound a disk $E^{'}$ in $A$. Thus
$A^{'}=(A-E)\cup E^{'}$ is also an annulus,  but $A^{'}$ can be
isotoped so that $|A^{'}\cap (D\cap V)|<|A\cap (D\cap V)|$, a
contradiction. Q.E.D. (Claim 3)

By Claim 3, there is an annulus $A$ which intersects the annulus
$D\cap V$ in only one arc. We may assume that $\partial D\subset
F\subset\partial_{-} V$. Now Let $N=N(A\cup(D\cap V))$ and $A^{*}$
be the closure of $\partial N(A\cup(D\cap V))-\partial_{-}
V\cup\partial_{+} V$. Then $A^{*}$ is also an annulus in $V$. We
may assume that $\partial_{1} A^{*}\subset
\partial_{+} V$ and $\partial_{2}A^{*}\subset F$. Since the genus
of $M=W\cup_{S} V$ is at least two, $\partial_{1} A^{*}$ is an
essential separating simple closed curve in $\partial_{+} V$ which
bounds a disk $B$ in $W$. Now there are two subcases:

{\bf Case 3.1.} \ $F$ is a torus.

In this case, $\partial_{2} A^{*}$ bounds a disk $B^{*}$ in $V$.
Now let $P=B\cup A^{*}\cup B^{*}$. Then $P$ is a 2-sphere which
intersects $\partial_{+} V$ in an essential simple closed curve.
That means that $M=W\cup V$ is a connected sum of two Heegaard
splittings.

By Lemma 3.1, $W\cup V$ is stabilized if and only if one of the
two factors is stabilized.

{\bf Case 3.2.} \ $g(F)\geq 2$.

Now $\partial_{2} A^{*}$ is an essential separating simple closed
curve in $\partial_{-} V$. $A^{*}\cup B$ is an essential disk
which intersects $\partial_{+} V$ in an essential simple closed
curve. By Case 1 and Case 2, $W\cup V$ is a boundary connected sum
 or a self-boundary connected sum of  Heegaard splittings .

Now by induction on $g(\partial_{-}V)$ and $g(\partial_{+} V)$, we
can prove Theorem 1. Q.E.D\vskip 3mm

{\bf Corollary 2.} Any unstabilized Heegaard splitting of a
handlebody $H$ is trivial.

{\bf Proof.} We shall proof this corollary by induction on the
genus of $H$. By [W], the unstabilized Heegaard splitting of a
3-ball is trivial.

Now suppose that $H$ is a handlebody of genus at least 1, and
$H=W\cup_F V$ is a unstabilized Heegaard splitting of $H$ such
that $W$ is a handlebody, $\partial H=
\partial_-V$ and  $g(F)>g(\partial_{-} V)$. Since $H$ is irreducible and
$\partial$-reducible, by Lemma 3.4, there is an essential disk $D$
such that $D$ intersects $F$ in an essential simple closed curve
in $F$. We may assume that $D\cap W$ is a disk and $D\cap V$ is an
annulus. Hence $\partial D\subset \partial_{-} V$.

If the genus of $\partial_-V$ is one, then by Case 3.1 in the
proof of Theorem 1, $H=W\cup_F V$ is a connected sum of two
Heegaard splittings, but $H$ is irreducible and $g(F)>1$, by [W],
$W\cup V$ is stabilized, a contradiction. So it must be that
$g(F)=g(\partial_- V)$, and the Heegaard splitting is trivial.

If the genus of $\partial_-V$ is larger then one, by Cases 1 and 2
in the proof of Theorem 1,  $W\cup V$ is stabilized, a
contradiction. So it must be that $g(F)=g(\partial_- V)$, and the
Heegaard splitting is trivial.
 Q.E.D.

Another proof of Corollary 2 is given by Fengchun Lei from
Scharlemann-Thompson's results.

 {\bf References.}\vskip 3mm

[CG]A. Casson and C. Gordon, \  Reducing Heegaard Splittings,
Topology and its Applications 27(1987) 275-283.

[H]W. Haken, \ Some results on surfaces in 3-manifolds, Studies in
Modern Topology (Math. Assoc. Amer., distributed by:
prentice-Hall, 1968), 34-98.

[L]F. Lei, \ On stability of Heegaard splittings, Math. Proc.
Cambridge Philos. Soc., 129(2000), 55-57.

 [Q]R. Qiu, \ Stabilization of Reducible Heegaard Splittings,
 Preprint.

[ST] M. Scharlemann and  A. Thompson, \ Heegaard splittings of
(surface) x I are standard, Math. Ann. 295 (1993) 549-564.

[W] F.Waldhausen, \ Heegaard-Zerlegungen der 3-sphere, Topology 7
(1968)195-203.

Jiming Ma

Department of Mathematics

Jilin University

Changchun, 130023

China.

E-mail:majiming2000@sohu.com \vskip 3mm

Ruifeng Qiu

Department of Mathematics

Dalian University of Technology

Dalian, 116024

China.

 qiurf@dlut.edu.cn\vskip 3mm
\enddocument